\documentclass[12pt,reqno]{amsart}

\newcommand{\mysection}[1]{\section{#1}
      \setcounter{equation}{0}}

\newtheorem{theorem}{Theorem}[section]
\newtheorem{lemma}[theorem]{Lemma}

\theoremstyle{definition}
\newtheorem{assumption}{Assumption}[section]

\theoremstyle{remark}

\newcommand\bL{\mathbb{L}}
\newcommand\bR{\mathbb{R}}

\newcommand\bW{\mathbb{W}}

\newcommand\cbrk{\hbox{\rm ]\kern-.15em]}}  
\newcommand\opar{\hbox{\rm \raise.1667ex\hbox{${\scriptstyle |
}$}\kern-.35em(}}

\newcommand\cB{\mathcal{B}}

\newcommand\cF{\mathcal{F}}

\newcommand\cP{\mathcal{P}}

\newcommand\cW{\mathcal{W}}

\begin{document}

\title[]
{SPDEs in divergence form with VMO
coefficients and filtering theory of partially observable
diffusion processes with Lipschitz coefficients}

\author[N.V.  Krylov]{N.V. Krylov}%
\thanks{This work  was partially supported
by NSF grant DMS-0653121}
\address{127 Vincent Hall, University of Minnesota,
Minneapolis,
       MN, 55455, USA}
\email{krylov@math.umn.edu}
 \keywords{Stochastic partial differential equations,
divergence equations, filtering equations}

\subjclass[2000]{60H15, 35R60}

\begin{abstract}
We present several results on  the smoothness in $L_{p}$ sense
of filtering densities under the Lipschitz continuity  assumption  
on the coefficients of a partially observable diffusion processes.
We obtain them by rewriting in divergence form
  filtering equation
which are usually considered in terms of formally adjoint to
operators  in nondivergence form.
 
\end{abstract}

\maketitle

\mysection{Introduction}

 Let $(\Omega,\cF,P)$ be a complete probability space
with an increasing filtration $\{\cF_{t},t\geq0\}$
of complete with respect to $(\cF,P)$ $\sigma$-fields
$\cF_{t}\subset\cF$. Denote by $\cP$ the predictable
$\sigma$-field in $\Omega\times(0,\infty)$
associated with $\{\cF_{t}\}$. Let
 $w^{k}_{t}$, $k=1,2,...$, be independent one-dimensional
Wiener processes with respect to $\{\cF_{t}\}$.  

We fix a stopping time $\tau$ and for $t\leq\tau$
in the Euclidean $d$-dimensional space $\bR^{d}$
of points $x=(x^{1},...,x^{d})$ we are considering the following
equation
\begin{equation}
                                          \label{11.13.1}
du_{t}=(L_{t}u_{t} 
+D_{i}f^{i}_{t}+f^{0}_{t})\,dt
+(\Lambda^{k}_{t}u_{t}+g^{k}_{t})\,dw^{k}_{t},
\end{equation}
  where $u_{t}=u_{t}(x)=u_{t}(\omega,x)$ is an unknown function,
$$
L_{t}\psi(x)=D_{j}\big(a^{ij}_{t}(x)D_{i}\psi(x)
+a^{j}_{t}(x)\psi(x)
\big)+b^{i}_{t}(x)D_{i}\psi(x)+c_{t}(x)\psi(x),
$$
$$
\Lambda^{k}_{t}\psi(x)=\sigma^{ik}_{t}(x)D_{i}\psi(x)
+\nu^{k}_{t}(x)\psi(x),
$$
the summation  convention
with respect to $i,j=1,...,d$ and $k=1,2,...$ is enforced
and detailed assumptions on the coefficients and the 
free terms will be given later.

One can rewrite \eqref{11.13.1} in the nondivergence form
assuming that the coefficients $a^{ij}_{t}$ and $a^{j}_{t}$
are differentiable
 in $x$ and then one could apply the results
from \cite{Kr99}. It turns out that the 
differentiability of $a^{ij}_{t}$ and $a^{j}_{t}$ is not needed
for the corresponding 
counterparts of the results in \cite{Kr99} to be true
and showing this and generalizing the 
corresponding results of \cite{Ki}
is one of the main purposes of   Section
\ref{section 7.10.1} of the present article.
We assume, roughly speaking, that $a^{ij}_{t}(x)$ are
measurable in $t$ and of class VMO with respect to $x$.

One of the main motivations for developing the theory  of SPDEs
comes from filtering theory of partially observable diffusion processes.This problem   is stated as follows.
 Let $d\geq1$, $d_{1}>d$ 
be integers .

Consider a $d_{1}$-dimensional two component process
 $z_{t}=(x_{t},y_{t})$
with $x_{t}$ being $d$-dimensional and $y_{t}$ 
$(d_{1}-d)$-dimensional. We
assume that $z_{t}$ is a diffusion process
 defined as a solution of the system
\begin{equation}\begin{split}         
                                            \label{eq3.2.14} 
& dx_{t}=b(t,z_{t}) dt+\theta (t,z_{t})dw_{t}, \\ 
& dy_{t}=B(t,z_{t}) dt+\Theta (t,y_{t})dw_{t}
\end{split}\end{equation}
with some initial data.

The coefficients of \eqref{eq3.2.14} are assumed to be
 vector- or matrix-valued
functions of appropriate dimensions defined on
 $[0,T]\times\bR^{d_{1}}$.
Actually $\Theta(t,y)$ is assumed to be independent
 of $x$, so  that it
is a function on $[0,T]\times\bR^{d_{1}-d}$ rather than 
$[0,T]\times\bR^{d_{1}}$
but as always we may think of $\Theta(t,y)$ as a
 function of $(t,z)$ as well.

The component $x_{t}$ is treated as unobservable
and $y_{t}$ as the only observations available.
The problem is to find a way to compute the density
$\pi_{t}(x)$
of the conditional 
distribution of $x_{t}$ given $y_{s},s\leq t$.
Finding an equation satisfied by $\pi_{t}$
(filtering equation) is considered
to be a solution of the (filtering) problem.
The filtering equations turn out to be 
particular cases of SPDEs.

In 1964 in \cite{Ku64} the filtering equations were proposed
in a somewhat nonrigorous way and  most likely some terms
in these equations
appeared from stochastic integrals written in
 the Stratonovich form and the others appeared from
  the It\^o integrals. Perhaps, the author
of \cite{Ku64} realized this too and published
an attempt to rescue some results of \cite{Ku64}
 in 1967 in \cite{Ku67}. This attempt
  turned  successful for  simplified models
  without the so-called cross terms.

Meanwhile, in 1966 in \cite{Sh66} the correct filtering
 equations in full generality,
yet assuming some regularity of the filtering density,
 were presented.
This is the reason we propose to call the filtering equations
in the case of partially observable diffusion
processes {\em Shiryaev's equations\/}
and their particular case without cross terms
{\em Kushner's equations\/}.

In case $d=1$ the result of \cite{Sh66} is presented
 in \cite{LS}
on the basis of  
the famous Fujisaki-Kallianpur-Kunita theorem (see \cite{FKK})
about the filtering equations in a very general setting.
 Some authors
even call the filtering equation for diffusion processes
the Fujisaki-Kallianpur-Kunita equation.

By adding to the Fujisaki-Kallianpur-Kunita theorem
 some simple facts from the theory of SPDEs,
the a priori regularity assumption  was removed 
  in \cite{KR} and under the Lipschitz  
and uniform nondegeneracy assumption
the $L_{2}$-version
of Theorem \ref{thm3.2.22} was proved. 
The basic result of \cite{KR} is that $\pi_{t}\in W^{1}_{2}$.
It is also proved that
if the coefficients are smoother, $\pi_{t}(x)$ is smoother too.
The nondegeneracy assumption was later removed 
(see \cite{Ro}) on the account of assuming
 that $\theta\theta^{*}$ is three
times continuously differentiable in $x$. It is again
proved that $\pi_{t}\in W^{1}_{2}$ and $\pi_{t}$ is even smoother
if the coefficients are smoother.

In \cite{Kr99} the results of \cite{KR} were improved,
$\theta\theta^{*}$ is assumed to be twice
  continuously differentiable in $x$ and it is shown that
$\pi_{t}\in W^{2}_{p}$ with any $p\geq2$.
 
The above mentioned results of \cite{KR}, \cite{Ro},
and \cite{Kr99}
use  filtering theory
in combination with the theory of SPDEs, the latter
being  stimulated by certain needs of  filtering theory. 
It turns out  that the theory of SPDEs alone can be used
to obtain the above mentioned regularity
 results about $\pi_{t}$ without knowing anything from 
 filtering theory itself. It also can be used to solve
other problems from  filtering theory.

The first ``direct'' (only using the theory of SPDEs)
proof of regularity of $\pi_{t}$
is given in \cite{KZ} in the case that
  system
\eqref{eq3.2.14} defines a nondegenerate diffusion process
and $\theta\theta^{*}$ is twice
 continuously differentiable in $x$. 
It is proved that $\pi_{t}\in W^{2}_{p}$ with any $p\geq2$
as in \cite{Kr99}.
Advantages of having arbitrary $p$ are seen from results
like our Theorem \ref{theorem 1.19.1}.
Of course, on the way of investigating $\pi_{t}$
in \cite{KZ} the filtering equations are derived 
``directly'' in an absolutely
different manner than before (on the basis of an idea
from \cite{KR81}).

In Section \ref{section 2.2.1} of this article we relax the smoothness assumption 
in \cite{KZ} to the assumption that the coefficients
of \eqref{eq3.2.14} are merely Lipschitz continuous,
the assumption which is almost always supposed to hold
when one deals with systems like \eqref{eq3.2.14}.
We find that $\pi_{t}\in W^{1}_{p}$.
Thus, under the weakest smoothness assumptions we obtain
the best (in the author's opinion)
 regularity result on $\pi_{t}$. 
In particular, we prove that if the initial data
is sufficiently regular, then the filtering density
is almost Lipschitz continuous in $x$ and $1/2$
H\"older continuous in $t$.
However, we still
assume $z_{t}$ to be nondegenerate.
Our approach is heavily based on analytic results.
There is also a probabilistic approach developed in \cite{Ku2}
and based on explicit formulas for solutions introduced in \cite{Pa}
and later developed
in \cite{KR81}  and \cite{Ku1} (also see references therein). This approach cannot give
as sharp results as ours in our situation.

It seems to the author that 
under the same assumptions of Lipschitz continuity,
 by following an idea
from 
  \cite{K79} one can solve another problem from filtering theory,
  the so-called innovation problem, and obtain the equality
$$\sigma\{y_{s},s\leq t\}
=\sigma\{\check{w}_{s},s\leq t\},
$$
where $\check{w}_{t}$ is the innovation Wiener process
of the problem (its definition is reminded
in Section \ref{section 2.2.1}). Recall that
for degenerate diffusion processes the positive solution
of the innovation problem is obtained in \cite{Pu}
again on the basis of the theory of SPDEs under the assumption
that the coefficients are more regular.

By the way, in our situation, if the coefficients are more regular,
the filtering
equation can be rewritten in a nondivergence form and then
additional smoothness of the filtering density,
existence of which is already established in this article,
is obtained on the basis of regularity results from \cite{Kr99}.

Although for the proof of the above 
mentioned
results concerning the filtering equations
it suffices to
use article \cite{Ki} about SPDEs in divergence form
 with continuous coefficients, we prefer to give  
 more general results
 borrowed from \cite{Div} in Section \ref{section 7.10.1}.
 In Section \ref{section 2.2.1} we present some results about
 the filtering equations from \cite{Fil}.
 
 We finish this section by introducing some notation.
 Let $K, 
\delta>0$ be fixed finite constants, $p\in[2,\infty)$.
 Denote $L_{p}=L_{p}(\bR^{d})$,
$C^{\infty}_{0}=C^{\infty}_{0}(\bR^{d})$. Introduce
$$
D_{i}=\frac{\partial}{\partial x^{i}},\quad i=1,...,d.
$$
By $Du$   we mean the gradient  with respect
to $x$ of a function $u$ on $\bR^{d}$.
As usual,  
$$
W^{1}_{p}=\{u\in L_{p}: Du\in L_{p}\},
\quad
 \|u\|_{W^{1}_{p}}=
\|u\|_{L_{p}}+\|Du\|_{L_{p}}.
$$

We use the same notation $L_{p}$ for vector- and matrix-valued
or else
$\ell_{2}$-valued functions such as
$g_{t}=(g^{k}_{t})$ in \eqref{11.13.1}. For instance,
if $u(x)=(u^{1}(x),u^{2}(x),...)$ is 
an $\ell_{2}$-valued measurable function on $\bR^{d}$, then
$$
\|u\|^{p}_{L_{p}}=\int_{\bR^{d}}|u(x)|_{\ell_{2}}^{p}
\,dx
=\int_{\bR^{d}}\big(
\sum_{k=1}^{\infty}|u^{k}(x)|^{2}\big)^{p/2}
\,dx.
$$

Recall that $\tau$ is a stopping time and introduce
$$
\bL _{p}(\tau):=L_{p}(\opar 0,\tau\cbrk,\cP,
L_{p}),\quad
\bW^{1}_{p}(\tau):=L_{p}(\opar 0,\tau\cbrk,\cP,
W^{1}_{p}).
$$
We also need the space $\cW^{1}_{p}(\tau)$,
which is the space of functions $u_{t}
=u_{t}(\omega,\cdot)$ on $\{(\omega,t):
0\leq t\leq\tau,t<\infty\}$ with values
in the space of generalized functions on $\bR^{d}$
and having the following properties:

(i)   $u_{0}\in L_{p}(\Omega,\cF_{0},L_{p})$;

(ii)    $u
\in \bW^{1}_{p}(\tau )$;

(iii) There exist   $f^{i}\in \bL_{p}(\tau)$,
$i=0,...,d$, and $g=(g^{1},g^{2},...)\in \bL_{p}(\tau)$
such that
 for any $\varphi\in C^{\infty}_{0}$ with probability 1
for all   $t\in[0,\infty)$
we have
$$
(u_{t\wedge\tau},\varphi)=(u_{0},\varphi)
+\sum_{k=1}^{\infty}\int_{0}^{t}I_{s\leq\tau}
(g^{k}_{s},\varphi)\,dw^{k}_{s}
$$
\begin{equation}
                                                 \label{1.2.1}
+
\int_{0}^{t}I_{s\leq\tau}\big((f^{0}_{s},\varphi)-(f^{i}_{s},D_{i}\varphi)
 \big)\,ds,
\end{equation}
where by $(f ,\varphi)$ we mean the action
of a generalized function $f$ on $\varphi$, in particular,
if $f$ is a locally summable,
$$
(f,\varphi)=\int_{\bR^{d}}f(x)\varphi(x)\,dx.
$$
Observe that, for any $\phi\in C^{\infty}_{0}$, the process
$(u_{t\wedge\tau},\phi)$ is $\cF_{t}$-adapted and (a.s.) continuous.

The reader can find 
in \cite{Kr99} a discussion of (ii) and (iii),
in particular, the fact that the series in \eqref{1.2.1}
converges uniformly in probability on every finite
subinterval of $[0,\tau]$.
In case that property (iii) holds, we write
\begin{equation}
                                       \label{12.3.1}
du_{t}=(D_{i}f^{i}_{t}+f^{0}_{t})\,dt
+g^{k}_{t}\,dw^{k}_{t}
\end{equation}
for $t\leq\tau$
and this explains the sense in which equation
\eqref{11.13.1} is understood. Of course, we still need to
specify appropriate assumptions on the coefficients
and the free terms in \eqref{11.13.1}.

 The work was partially supported by NSF Grant DMS-0653121.
 
 \mysection{SPDEs in divergence form with VMO coefficients}
                                                                                                  \label{section 7.10.1} 

We are considering  \eqref{11.13.1}  under the following assumptions.
 
\begin{assumption}
                                        \label{assumption 1.2.1}

(i) The coefficients $a^{ij}_{t}$, $a^{i}_{t}$, $b^{i}_{t}$,
$\sigma^{ik}_{t}$, $c_{t}$, and $\nu^{k}_{t}$ are measurable 
with respect to $\cP\times \cB(\bR^{d})$, where $\cB(\bR^{d})$
is the Borel $\sigma$-field on $\bR^{d}$.

(ii) For all values of indices and arguments
$$
 |a^{i}_{t}|+|b^{i}_{t}|+|c_{t}|+|\nu|_{\ell_{2}}\leq
K,\quad c_{t}\leq0.
$$

(iii) For all values of the arguments and $\xi\in\bR^{d}$
\begin{equation}
                                             \label{1.3.2}
 a^{ij}_{t}  
  \xi^{i}
\xi^{j}\leq\delta^{-1}|\xi|^{2},\quad
(a^{ij}_{t}  
-  \alpha^{ij}_{t}) \xi^{i}
\xi^{j}\geq\delta|\xi|^{2},
\end{equation}
where   $\alpha^{ij}_{t}=(1/2)(\sigma^{i\cdot},\sigma^{j\cdot})
_{\ell_{2}}$.  
\end{assumption} 

It is worth emphasizing that we do not require
the matrix $(a^{ij})$ to be symmetric.
Assumption \ref{assumption 1.2.1} (i) guarantees
that equation \eqref{11.13.1} makes perfect sense
if $u\in\cW^{1}_{p}(\tau)$.

For   functions $h_{t}(x)$ on $[0,\infty)\times\bR^{d }$
and balls $B$ in $\bR^{d}$ introduce
$$
 h_{ t(B )} 
=\frac{1}{|B |}\int_{B }h_{t}( x)\,dx,
$$
where   $|B|$ is the volume of $B$. If
$\rho\geq0$, set $B_{\rho}=\{x:|x|<\rho\}$ and for
locally integrable $h_{t}(x)$ and   continuous
 $\bR^{d}$-valued function $x_{r},r\geq0$, introduce
$$
\text{osc}_{\rho} \, (h,x_{\cdot})=
\sup_{s\geq0}\frac{1}{\rho^{2}}
\int_{s}^{s+\rho^{2}}(|h_{r} -h_{ r(B+x_{r})}|)_{(B+x_{r})} \,dr,
$$  
where $B=B_{\rho}$.  Also for $y\in\bR^{d}$ set
$$
\text{Osc}_{\rho} \, (h,y)=
\sup_{|x_{\cdot}|_{C}\leq\rho} \sup_{r\leq\rho}\text{osc}_{r} 
\, (h,y+x_{\cdot}),
$$
where $|x_{\cdot}|_{C}$ is the sup norm of $|x_{\cdot}|$.
 Observe that $\text{ocs}_{\varepsilon}h=0 $ if $h_{t}(x)$ is
independent of $x$.

Denote by  $\beta_{0}$   one third
of the constant $\beta_{0}(d,p,\delta)>0 $
 from Lemma 5.1 of~\cite{Div}.

 \begin{assumption}
                                       \label{assumption 1.2.6}
                                       There exist a constant  $\varepsilon
\in(0,1]$ 
such that for any   $y\in\bR^{d}_{+}$
(and $\omega$) 
we have 
\begin{equation}
                                                   \label{8.7.1}
\text{Osc}_{\varepsilon } \,(a^{ij},y)\leq \beta_{0},
\quad\forall i,j.
\end{equation}

Furthermore,
$$
(a^{jk}_{t}(x)  
-  \alpha^{jk}_{t}(y)) \xi^{j}\xi^{k}\geq\delta|\xi|^{2}
$$
for all   $t $, $\xi$, and $x$  satisfying
  $|x-y|\leq\varepsilon $.

\end{assumption}

Let $\beta_{1}=\beta_{1}(d,p,\delta,\varepsilon )>0$
be the constant   
from  Lemma 5.2 of \cite{Div}.

\begin{assumption} 
                                       \label{assumption 8.10.1}
There exists a constant  $\varepsilon_{1}>0$ 
such that for any $t\geq0$  
we have
$$
|\sigma^{i\cdot}_{t}(x)-\sigma^{i\cdot}_{t}(y)|_{\ell_{2}}\leq \beta_{1},
$$
whenever 
$x,y\in\bR^{d}_{+}$,  $ |x-y|\leq \varepsilon_{1}$, $ i=1,...,d$ .

\end{assumption} 

Finally, we describe the space of initial data.
Recall that for $p\geq2$
the Slobodetskii space $W^{1-2/p}_{p}
=W^{1-2/p}_{p}(\bR^{d})$ of functions $u_{0}(x)$
can be introduced
as the space of traces on $t=0$ of 
(deterministic) functions $u $
such that
$$
u\in L_{p}(\bR_{+},W^{1}_{p}),\quad\partial u/\partial t
\in L_{p}(\bR_{+},H^{-1}_{p}),
$$
where $\bR_{+}=(0,\infty)$ and $H^{-1}_{p}=(1-\Delta)^{-1/2}L_{p}$.
For such functions there is a 
(unique) modification denoted again
$u$ such that $u_{t}$ is a continuous $L_{p}$-valued
function on $[0,\infty)$ so that $u_{0}$ is well defined.
Any such $u_{t}$ is called an extension of $u_{0}$.

The norm in $W^{1-2/p}_{p}$ can be defined as
the infimum  of
$$
\|u\|_{ L_{p}(\bR_{+},W^{1}_{p})}
+\|\partial u/\partial t \|_{L_{p}(\bR_{+},H^{-1}_{p})}
$$
over all extensions $u_{t} $ of elements $u_{0} $.

\begin{theorem}
                                    \label{theorem 12.7.1} 
Let    
$f^{j},g\in\bL_{p}(\tau)$  and let
$u_{0}\in L_{p}(\Omega,\cF_{0},
 W^{1-2/p}_{p})$.   Then

(i) 
Equation \eqref{11.13.1}
for $t\leq T\wedge\tau$
has a unique solution  
$u\in\cW^{1}_{p}(T\wedge\tau)$ with initial data $u_{0}$
for any constant $T\in(0,\infty)$.
 
(ii) There exists a set $\Omega'\subset
\Omega$ of full probability
such that $u_{t\wedge\tau}I_{\Omega'}$ is a continuous 
$\cF_{t}$-adapted 
$L_{p}$-valued functions
of $t\in[0,\infty)$.
 
\end{theorem}

 Assertion (ii) of Theorem \ref{theorem 12.7.1}  follows from assertion (i)
 and Theorem \ref{theorem 12.3.1}.

 Here is a result about continuous dependence
of solutions on the data.

\begin{theorem}
                                    \label{theorem 1.4.1}
 
Assume that for each $n=1,2,...$
we are given   functions $a^{ij}_{nt}$, $a^{i}_{nt}$, $
b^{i}_{nt}$, $c_{nt}$, $\sigma^{ik}_{nt}$, $\nu^{k}_{nt}$, 
$f^{j}_{nt}$, $g^{k}_{nt}$, and $u_{n0} $
having the same meaning  and satisfying the same assumptions
with the same $\delta,K$, $\varepsilon$, $\varepsilon_{1}$,  $
\beta_{0}$, and $\beta_{1}$ as the original ones.
Assume that for 
$i,j=1,...,d$ and
almost all $(\omega,t,x)$ we have
$$
(a^{ij}_{nt},a^{i}_{nt},b^{i}_{nt},c_{nt})\to
(a^{ ij}_{t},a^{ i}_{t},b^{ i}_{t},c _{t}),
$$
$$
|\sigma^{ i\cdot}_{nt}-\sigma^{ i\cdot}_{t}|_{\ell_{2}}+
|\nu _{nt}-\nu _{t}|_{\ell_{2}}\to0, 
$$
as $n\to\infty$. Also  assume that 
$$
\sum_{j=0}^{d}(\|f^{ j}_{n}-f^{j}\|_{\bL_{p}(\tau)}
+
\|g_{n }-g\|_{\bL_{p}(\tau)}
 +\|u_{n0} -u_{0}\|_{L_{p}(\Omega,\cF_{0},
 W^{1-2/p}_{p})}  \to0
$$
as $n\to\infty$. Let $u_{n}$
 be the unique solutions
of equations \eqref{11.13.1} for $t\leq\tau$ constructed from
$a^{ij}_{nt}$, $a^{i}_{nt}$, $
b^{i}_{nt}$, $c_{nt}$, $\sigma^{ik}_{nt}$, $\nu^{k}_{nt}$, 
$f^{j}_{nt}$, and $g^{k}_{nt}$ and having initial
values  $u_{n0}  $. 

Then,  for any $T\in[0,\infty)$ as $n
\to\infty$, we have
$\|u_{n}-u\|_{\bW^{1}_{p}(T\wedge\tau )}\to0$ and
$$
E\sup_{t\leq\tau\wedge T}
\|u_{n t}-u_{t}\|_{L_{p}}^{p}\to0.
$$
 \end{theorem}

 In many situation the following maximum principle
 based on the results of \cite{Kr07}
is useful.
 \begin{theorem}
                                    \label{theorem 1.4.3}

Suppose that,  for $q\in[2,p]$, Assumptions \ref{assumption 1.2.6}
and \ref{assumption 8.10.1} are satisfied with $\beta_{0} \leq\beta_{0}(d,q,\delta)$ 
and $\beta_{1}\leq\beta_{1}(d,q,\delta,\varepsilon)$.
Also suppose that  
$u_{0}\in L_{p}(\Omega,\cF_{0},
 W^{1-2/q}_{q})$, $q\in[2,p]$,
$u _{0}\geq0$,
$f^{i}=0$, $i=1,...,d$, $f^{0}\geq0$, $g=0$.
Then for the solution $u$ almost surely
we have $u_{t}\geq0$ for all finite $t\leq\tau$.
 
\end{theorem}

Part of the proofs of the above results is based on the following
It\^o's formula.

\begin{theorem}
                                     \label{theorem 12.3.1}
  Let
$u\in\cW^{1}_{p}(\tau)$, 
$f^{j}\in\bL_{p}(\tau)$,
  $g=(g^{k})\in\bL_{p}(\tau)$ and assume that
\eqref{12.3.1} holds
for $t\leq\tau$ in the sense of generalized functions.
Then there is a set $\Omega'\subset\Omega$ of full probability
such that

(i) $u_{t\wedge\tau}I_{\Omega'}$
is a continuous $L_{p}$-valued $\cF_{t}$-adapted function on
$[0,\infty)$;

(ii) for all $t\in[0,\infty)$ and $\omega\in\Omega'$ It\^o's formula holds:
$$
\int_{\bR^{d}}|u_{t\wedge\tau}|^{p}\,dx
=\int_{\bR^{d}}|u_{0}|^{p}\,dx
+p
\int_{0}^{t\wedge\tau }\int_{\bR^{d}}|u _{s}|^{p-2}
u _{s}
g^{k }_{s}\,dx\,dw^{k}_{s}
$$
$$+
\int_{0}^{t\wedge\tau }
\big( \int_{\bR^{d}}\big[p|u_{t}|^{p-2}u_{t}f^{0}_{t}
-p(p-1)|u_{t}|^{p-2}f^{i}_{t}D_{i}u_{t}
$$
\begin{equation}
                                            \label{4.19.5}
+(1/2)p(p-1)|u_{t}|^{p-2}|g_{t}|_{\ell_{2}}^{2}
\big]\,dx\big)\,dt.
\end{equation}

Furthermore,  for any $T\in[0,\infty)$ 
$$
 E\sup_{t\leq\tau\wedge T}
\|u_{t}\|^{p}_{L_{p}}\leq  2E\|u_{0}\|^{p}_{L_{p}}+
NT^{p-1}\|f^{0}\|^{p}_{\bL_{p}(\tau)}
$$
\begin{equation}
                                         \label{4.11.5}
+NT^{(p-2)/2}(\sum_{i=1}^{d}\|f^{i}\|^{p}_{\bL_{p}(\tau)}
+\|g\|^{p}_{\bL_{p}(\tau)}+\|Du\|^{p}_{\bL_{p}(\tau)}) ,
\end{equation}
where $N=N(d,p)$.  

\end{theorem}
 
 We have a direct proof of this result. However, \eqref{4.19.5}
 can also  be obtained by extending some arguments from \cite{Br}.

\mysection{Filtering equations}
                                         \label{section 2.2.1}
 Fix a constant $T\in(0,\infty)$
   and                                      
for simplicity  assume that
 $w_{t}$ in \eqref{eq3.2.14}
 is finite dimensional.
First we state and discuss our assumptions.

\begin{assumption}
                                             \label{asm3.2.15}
The functions $b$, $\theta $, $B$, and $\Theta $ are
 Borel measurable and bounded
functions of their arguments. Each of them satisfies 
the Lipschitz 
condition in $z$ with the constant $K $. 
\end{assumption}

\begin{assumption}
                                      \label{assumption 1.15.1} 
The process $z_{t}$ is uniformly nondegenerate:
for any $\lambda,z\in\bR^{d_{1}}$ and $t\in[0,T]$ we have
$$
 \tilde{a}^{ij}_{t}(z)\lambda^{i}
\lambda^{j}\geq\delta|\lambda|^{2}, 
$$
where $ 2\tilde{a} _{t}(z)=2(\tilde{a}^{ij}_{t}(z))=\theta(t,z)\theta^{*}(t,z)
+\Theta(t,y)\Theta^{*}(t,y)$.

\end{assumption}

Traditionally, Assumption \ref{assumption 1.15.1}
is split into two following assumptions 
the combination of which is equivalent to Assumption \ref{assumption 1.15.1} and
in which some
useful objects are introduced. These assumptions
 were also used in the past to reduce system \eqref{eq3.2.14}
to the so-called triangular form by replacing
$w_{t}$ with a different Brownian motion.

\begin{assumption}
                                           \label{asm3.2.16}
The symmetric matrix $\Theta \Theta^{*}$ is 
invertible and 
$$
\Psi :=(\Theta
\Theta^{*} )^{-\frac{1}{2}}
$$ 
is a bounded function of $(t,y)$.
\end{assumption}
 
\begin{assumption}\label{asm3.2.17}
For any $\xi \in \bR^{d}$, $z=(x,y)\in \bR^{d_{1}}$, and $t>0$,
we have
\[
|Q(t,y)\theta^{*}(t,z)\xi |^{2}\geq \delta |\xi |^{2},
\]
where $Q$ is the orthogonal projector on $\text{Ker}\,\Theta$. 
In other words,
\begin{equation}                        \label{eq3.2.17.1}
(\theta (I-\Theta^{*}\Psi^{2}\Theta )
\theta^{*} \xi ,\xi )\geq 
\delta |\xi|^{2}.
\end{equation}
\end{assumption}

\begin{assumption}                            \label{asm3.2.18}
The random vectors $x_{0}$ and $y_{0}$ are 
independent of the process
$w_{t}$. The conditional distribution of $x_{0}$ given $y_{0}$ has a
density, which we denote by $\pi_{0}(x)=\pi_{0}(\omega ,x)$. We have
$\pi_{0} \in L_{p}(\Omega , W_{p}^{1-2/p} )$.
\end{assumption}

Next we introduce few more notation.
Let
$$
\Psi_{t}=\Psi(t,y_{t}),\quad\Theta_{t}=\Theta(t,y_{t}),\quad
a_{t}(x) =\frac{1}{2}\theta \theta^{*}(t,x,y_{t}),
\quad b_{t}(x)=b(t,x,y_{t}),
$$
$$
\sigma_{t}(x) =\theta(t,x,y_{t}) \Theta^{*}_{t}\Psi_{t},\quad
\beta_{t}(x) =\Psi_{t}B(t,x,y_{t}).
$$
In the remainder of the article we use the notation
$$
D_{i}=\frac{\partial}{\partial x^{i}}
$$
only for $i=1,...,d$ and set
\begin{equation}                             \label{eq3.2.19.2}
L_{t}( x) = a^{ij}_{t}(x)D_{i}D_{j}+
 b^{i}_{t}(x)D_{i}\,,
\end{equation}
$$
L^{*}_{t}( x)u_{t}(x) = 
D_{i}D_{j}(
a^{ij}_{t}( x)u_{t}(x) )
- D_{i}(b^{i}_{t}( x)u_{t}(x) )
$$
\begin{equation}                     \label{e3.2.19.2}
=D_{j}\big(
a^{ij}_{t}(x) D_{i}u_{t}(x)
-b^{j}_{t}(x)u_{t}(x)+u_{t}(x)D_{i}
a^{ij}_{t}(x)\big),
\end{equation}
\begin{equation}                     \label{1.27.8}
\Lambda^{k }_{t}( x)u_{t}(x)  =\beta^{k}_{t}( x)u_{t}(x) +
 \sigma^{ik}_{t}( x) D_{i}u_{t}(x),
\end{equation}
$$
\Lambda^{k*}_{t}( x)u_{t}(x)  =\beta^{k}_{t}( x)u_{t}(x) -
D_{i}(\sigma^{ik}_{t}( x) u_{t}(x) )
$$
\begin{equation}\label{eq3.2.19.3}
=-\sigma^{ik}_{t}(x)D_{i}u_{t}(x)
+(\beta^{k}_{t}(x)-D_{i}
\sigma^{ik}_{t}(x))u_{t}(x),
\end{equation}
where $t\in[0,T]$,   $x\in \bR^{d}$, $k=1,...,d_{1}-d$,
and as above
 we use the summation convention.
Observe that Lipschitz continuous functions have
bounded generalized derivatives and by
$$
D_{i}
a^{ij}_{t},\quad D_{i}
\sigma^{ik}_{t}
$$
we mean these derivatives. Obviously, the operator
$L$ defined by \eqref{eq3.2.19.2} is uniformly elliptic 
with constant
of ellipticity $\delta$.

Finally, by $\cF_{t}^{y}$ we denote the completion of 
$\sigma \{ y_{s}:s\leq t\}$
with respect to $P,\cF$.

Let us consider the following initial value problem
\begin{equation}
                                        \label{eq3.2.20}
d\bar{\pi }_{t}(x)=L^{*}_{t}(x )\bar{\pi }_{t}(x)\,dt+
 \Lambda^{k*}_{t}(x )\bar{\pi }_{t}(x)
 \Psi^{kr}_{t}\,dy^{r}_{t},
\end{equation}
$$
\bar{\pi }_{0}(x)= \pi  _{0}(x),
$$
where $t\in[0,T]$, $x\in\bR^{d}$, and
$\bar{\pi }_{t}(x)=\bar{\pi }_{t}(\omega,x)$. 
Equation \eqref{eq3.2.20}
is called the Duncan-Mortensen-Zakai or just the Zakai equation.

We understand this equation and the initial
condition in the following sense.
We are looking for a function $\bar{\pi}=\bar{\pi}_{t}(x)
=\bar{\pi}_{t}(\omega,x)$, $\omega\in\Omega$,
$t\in[0,T]$, $x\in\bR^{d}$, such that

(i) For each $(\omega,t)$, $\bar{\pi}_{t}(\omega,x)$
is a generalized function on $\bR^{d}$,

(ii) We have $\bar{\pi}\in L_{p}(\Omega\times
[0,T],\cP,W^{1}_{p} )$,

(iii) For each $\varphi\in C^{\infty}_{0}(\bR^{d})$
with probability one for all $t\in[0,T]$ it holds that
$$
(\bar{\pi}_{t},\varphi)=(\pi_{0},\varphi)
-\int_{0}^{t}(
a^{ij}_{t} D_{i}\bar{\pi}_{t}
-b^{j}_{t}\bar{\pi}_{t}+\bar{\pi}_{t}D_{i}
a^{ij}_{t},D_{j}\varphi)\,dt
$$
\begin{equation}
                                                 \label{1.16.1}
-\int_{0}^{t}
(\sigma^{ik}_{t}D_{i}\bar{\pi}_{t}
+( D_{i}
\sigma^{ik}_{t}-\beta^{k}_{t})\bar{\pi}_{t},\varphi)\Psi^{kr}_{t}
\big(B^{r}(t,z_{t})\,dt+\Theta^{rs}(t,y_{t})\,dw^{s}_{t}\big).
\end{equation}

Observe that all expressions in \eqref{1.16.1} are well defined
due to the fact that the coefficients 
of $\bar{\pi}$ and of  
$D_{i}\bar{\pi}$
are bounded
and appropriately measurable and $\bar{\pi},
D_{i}\bar{\pi}\in L_{p}(\Omega\times[0,T],
\cP,L_{p} )$.

Hence, equation \eqref{eq3.2.20} has the same form as \eqref{11.13.1}
and  
the existence and uniqueness part of   Lemma 
\ref{lm3.2.21} below
  follow  from Theorem \ref{theorem 12.7.1}. The second assertion of the lemma
follows from Theorem \ref{theorem 1.4.3}.

\begin{lemma}
                                            \label{lm3.2.21} 
There exists a unique solution $\bar{\pi}$ 
of \eqref{eq3.2.20} with initial
condition $\pi_{0}$
in the sense explained above. In addition,   
$
\bar{\pi}_{t}\geq0
$
for all $t\in[0,T]$ (a.s.).
\end{lemma}

Here is a basic result of filtering theory
for partially observable diffusion
processes. Its relation to the previously
known ones is discussed above.

\begin{theorem}
                                               \label{thm3.2.22}
Let $\bar{\pi }$ be the function from Lemma \ref{lm3.2.21}.
Then 
\begin{equation} 
                                             \label{1.27.07} 
0<\int_{\bR^{d}}\bar{\pi }_{t}(x)\,dx=(\bar{\pi }_{t},1)<\infty 
\end{equation}
 for 
all $t\in[0,T]$  (a.s.)
and for any $t\in[0,T]$ and real-valued,
bounded or nonnegative, (Borel) measurable function $f$ 
given on $\bR^{d}$
\begin{equation} 
                                             \label{10.13.3}
E[f(x_{t})|\cF_{t}^{y}]=
\frac{(\bar{\pi }_{t},f)}
{(\bar{\pi }_{t},1)}\quad
\text{(a.s.).}
\end{equation}
\end{theorem}

Equation \eqref{10.13.3} shows (by definition) that
$$
\pi_{t}(x):=\frac{\bar{\pi }_{t}(x) }
{(\bar{\pi }_{t},1)}
$$
is a conditional density of distribution of
$x_{t}$ given $y_{s},s\leq t$. Since, generally,
$(\bar{\pi }_{t},1)\ne1$,  one calls $\bar{\pi }_{t}$
an unnormalized conditional density of distribution of
$x_{t}$ given $y_{s},s\leq t$.

We derive Theorem \ref{thm3.2.22}
from Theorem \ref{theorem 1.4.1} and the result of
\cite{KZ} where more regularity on the coefficients
is assumed.

The following is a direct corollary of 
embedding theorems from \cite{Kr99}.

\begin{theorem}
                                        \label{theorem 1.19.1}
Let $\pi_{0}$ be a nonrandom 
function and $\pi_{0}\in W^{1-2/p}_{p} $
for all $p\geq2$, which happens
for instance, if $\pi_{0}$ is a
Lipschitz
continuous function with compact support.
Then for any $\varepsilon\in(0,1/2)$ almost surely
$\bar{\pi}_{t}(x)$ is $1/2-\varepsilon$ H\"older 
continuous in $t$ with a constant independent of $x$,
$\bar{\pi}_{t}(x)$ is $1-\varepsilon$ H\"older 
continuous in $x$ with a constant independent of $t$,
and the above mentioned (random) constants have all moments.

\end{theorem}

 \bibliographystyle{plain}

\end{document}